\documentclass[12pt,reqno]{amsart}
\usepackage{graphicx}
\usepackage{epsfig}
\usepackage{amssymb,amsthm,amsmath,amsfonts,amstext,mathrsfs}
\usepackage{latexsym}
\usepackage{url}
\usepackage{palatino}
\newtheorem{thm}{Theorem}

\newtheorem{prop}{Proposition}
\newtheorem{conj}{Conjecture}

\input xy
\xyoption{all}

\def\d{\,{\rm{d}}}

\title[Formula for the moments of $?(x)$]
{Semi-regular continued fractions and an exact formula for the moments of the Minkowski question mark function}
\author[Giedrius Alkauskas]{ Giedrius Alkauskas}
\begin{document}
\begin{abstract} This paper continues investigations on various integral transforms of the Minkowski question mark function. In this work we finally establish the long-sought formula for the moments, which does not explicitly involve regular continued fractions, though it has a hidden nice interpretation in terms of semi-regular continued fractions. The proof is self-contained and does not rely on previous results by the author.
\end{abstract}
\maketitle
\begin{center}
\end{center}
\begin{center}
\end{center}
\section{Introduction and main result}
\footnotetext[1]{The author gratefully acknowledges support
from the Austrian Science Fund (FWF) under the project Nr. P20847-N18.}
This work is the fifth paper devoted to investigations of various integral transforms of the function in title. Other four are \cite{alkauskas1},
\cite{alkauskas2}, \cite{alkauskas3} and \cite{alkauskas4} in chronological order of research.\\

The Minkowski question mark function $?(x)$ is defined by
\begin{eqnarray}
?([0,a_{1},a_{2},a_{3},...])=2^{1-a_{1}}-2^{1-(a_{1}+a_{2})}+2^{1-(a_{1}+a_{2}+a_{3})}-...,\quad x\in[0,1],
\label{min}
\end{eqnarray}
where $x=[0,a_{1},a_{2},...]$ stands for the representation of $x$ by a regular continued fraction. In case $x$ is rational, we have two developments as a continued fraction, but they both, as it is easy to check, produce the same value for $?(x)$. This function has many fascinating properties; it is a deep conviction of the author that many of them are still waiting to be discovered.  In the past, this function was studied from many different standpoints: dynamical systems, multi-fractal analysis \cite{kessebohmer}, metric number theory, complex dynamics, analytic number theory, combinatorics. It can be generalized in many ways, including generalizations to higher dimensions \cite{panti}. From the point of view of number theory, maybe the most interesting generalization was given in \cite{grabner_tichy} and continued in \cite{lamberger}. We refer to the paper \cite{alkauskas1} for an overview of the existing literature. The internet page \cite{mink} contains an extensive bibliography list on the Minkowski question mark function.\\

The  function $?(x)$ is continuous, strictly increasing bijection of $[0,1]$, it is singular ($?'(x)=0$ almost everywhere), and  for $x\in[0,1]$ it satisfies the following functional equations:
\begin{eqnarray*}
?(x)=\left\{\begin{array}{c@{\qquad}l} 1-?(1-x),
\\ 2?(\frac{x}{x+1}). \end{array}\right.
\end{eqnarray*}
\indent The main idea of previous research by the author was that these functional equations carry onto the various integral transforms (Laplace, Stieltjes, Mellin) of the question mark function, and the so-obtained objects are unique (subject to certain natural regularity conditions), and thus they arise exactly from the Minkowski question mark function. Truly, the main motivation for the introduction of these integral transforms was to get a better understanding of the structure of \emph{moments} (which are analogues and share many properties with the \emph{periods} of Maass wave forms). For an integer $L\geq 0$, the moments are defined by
\begin{eqnarray*}
m_{L}=\int\limits_{0}^{1}x^{L}\d?(x).
\end{eqnarray*}
\indent The Minkowski question mark function describes the limit distribution of the rational numbers in the Farey tree. Each generation of the Farey tree consists of those rational numbers whose sum of the elements of the continued fraction is fixed. Thus, these moments can  also be defined by
\begin{eqnarray}
m_{L}=\lim_{n\rightarrow\infty}2^{2-n}\sum\limits_{a_{1}+\ldots+a_{s}=n\atop a_{i}\geq 1,\quad a_{s}\geq 2}[0,a_{1},\ldots,a_{s}]^{L}.
\label{riba}
\end{eqnarray}
This gives the second formula. There exist several seemingly unrelated procedures which also produce the moments. We will mention one of them, a conjectural and rather mysterious procedure (it assumes that a certain analytic function $G({\sf p},z)$ in two complex variables, whose set of irregularities for ${\sf p}$ fixed is a fractal Julia set, has an analytic continuation), and refer the reader to \cite{alkauskas4} for the further details.
\begin{conj}
Define the rational functions  $\mathbf{Q}_{n}(z)$, $n\geq 0$ (with rational coefficients, which are $p-$adic integers for $p\neq 2$) by
\begin{eqnarray*}
\mathbf{Q}_{0}(z)=-\frac{1}{2z},\text{ and recurrently by }
\mathbf{Q}_{n}(z)=\frac{1}{2}\sum\limits_{j=0}^{n-1}\frac{1}{j!}\cdot\frac{\d^{j}}{\d
z^{j}}\mathbf{Q}_{n-j-1}(-1)\cdot\Big{(}z^{j}-\frac{1}{z^{j+2}}\Big{)}.
\end{eqnarray*}
Then $\Lambda(t)=\sum\limits_{n=0}^{\infty}\frac{\mathbf{Q}^{'}_{n}(-1)}{n!}t^{n}$ is an entire function and $m_{2}=\int_{0}^{\infty}\Lambda(t)e^{-t}\d t$.
\end{conj}
The sequence $\mathbf{Q}^{'}_{n}(-1),n\in\mathbb{N}_{0}$, begins with $\frac{1}{2},-\frac{1}{2},1,-\frac{5}{2},\frac{25}{4},-16,43,-\frac{971}{8},\frac{1417}{4},\ldots$
Note that in this conjecture the continued fractions are also involved implicitly, since a shift by $-1$ in the formula which defines $\mathbf{Q}_{n}(z)$ means that we employ the function $z\rightarrow z-1$, and  we obviously also employ $z\rightarrow\frac{1}{z}$, which together are intricately related with continued fractions.\\

In this paper we finally present an exact formula for the moments, which does not explicitly involve continued fractions; rather, it has a nice hidden explanation in terms of semi-regular ones. Surprisingly, we can present a self-contained treatment which neither uses the functional equations for $?(x)$, nor is employing any previous results; though, of course, the generating function of moments (the so-called \emph{dyadic period function}, a close relative of period functions associated with Maass wave forms), the three-term functional equation it satisfies, and some other objects are implicitly present.
\begin{thm} For an integer $L\geq 1$, the moments of the Minkowski question mark function can be expressed by the following exact formula
\begin{eqnarray*}
\int\limits_{0}^{1}x^{L}\d ?(x)=\frac{1}{(L-1)!}\sum\limits_{\ell=0}^{\infty}\quad\int\limits_{[0,\infty)^{\ell+1}}x_{0}^{L}\cdot
\frac{(x_{0}x_{\ell})^{-1/2}\cdot\prod\limits_{i=0}^{\ell-1}I_{1}(2\sqrt{x_{i}x_{i+1}})}{\prod\limits_{i=0}^{\ell}e^{x_{i}}(2e^{x_{i}}-1)}\d x_{0}\ldots\d x_{\ell}.
\end{eqnarray*}
Here $I_{1}(x)=J_{1}(ix)$ stands for the Bessel $I-$function. The $\ell-$th term of this series is of magnitude $\ll 2^{-\ell}$.
\end{thm}
\noindent (By a convention, an empty product is assumed to be $1$. This applies to the case $\ell=0$). Unfortunately, this expression is not very suitable for calculations. The manipulations with matrices of high order can produce much better numerical values for these constants \cite{alkauskas4}. Further, the symmetry property of $?(x)$, given by $?(x)+?(1-x)=1$, produces linear relations among the moments. For example, $m_{1}=\frac{1}{2}$, $3m_{2}-2m_{3}=\frac{1}{2}$. Since the symmetry property is obvious if one is dealing with the regular continued fractions and does not manifest in semi-regular ones, even the fact that for $L=1$ the right-hand-side is $\frac{1}{2}$ is a strange by-product identity. Thus, the first four terms give (with $10$ precise decimal digits)
\begin{eqnarray*}
0.3862943611+0.0791502471+0.0226858500+0.0074990924=0.4956295506.
\end{eqnarray*}
Consequently, it is not excluded that there might be a simpler expression for the moments.
\section{The proof}
We will show that the derived formula has an interpretation in terms of \emph{semi-regular continued fractions}. These are defined by
\begin{eqnarray*}
[[b_{1},b_{2},b_{3},\ldots]]=\cfrac{1}{b_{1}-\cfrac{1}{b_{2}-\cfrac{1}{b_{3}-\ddots}}},
\end{eqnarray*}
where integers $b_{i}\geq 2$. Each real irrational number $x\in(0,1)$ has a unique representation in this form, and rationals $x\in(0,1)$ have two of them: one finite, one infinite which ends in $[[2,2,2,\ldots]]$.\\
\indent The following proposition establishes the relation between semi-regular continued fractions and $?(x)$, since this relation easily follows from the results obtained in \cite{ryde}; see also the concise exposition \cite{ramharter}.
\begin{prop} We have
\begin{eqnarray}
?\Big{(}[[b_{1},b_{2},b_{3},\ldots]]\Big{)}=2^{1-b_{1}}+2^{2-(b_{1}+b_{2})}+2^{3-(b_{1}+b_{2}+b_{3})}+\ldots\label{sum}
\end{eqnarray}
\end{prop}

{\it Proof.} As before, two developments for rational numbers produce the same value. Denote the function on the right by $v(x)$, $x=[[b_{1},b_{2},b_{3},\ldots]]$. Let $2_{a}$, $a\in\mathbb{N}_{0}$, denote the string $\underbrace{2,2,\ldots,2}_{a}$. Then, as was proved in \cite{ramharter}, we have:
\begin{eqnarray*}
x=[0,a_{1},a_{2},a_{3},\ldots]=[[a_{1}+1,2_{a_{2}-1},a_{3}+2,2_{a_{4}-1},a_{5}+2,2_{a_{6}-1},\ldots]].
\end{eqnarray*}
(Note that $a_{1}$ is slightly exceptional). Let us calculate $v(x)$ using the expansion (\ref{sum}).  The first $a_{2}$ terms of this expansion contribute
\begin{eqnarray*}
2^{-a_{1}}+\sum\limits_{i=1}^{a_{2}-1}2^{i+1-(a_{1}+1)-2i}=2^{-a_{1}}+2^{-a_{1}}(1-2^{1-a_{2}})=2^{1-a_{1}}-2^{1-(a_{1}+a_{2})}.
\end{eqnarray*}
Further, the direct calculation shows that the first $a_{2}+a_{4}$ terms of the expansion (\ref{sum}) produce
$2^{1-a_{1}}-2^{1-(a_{1}+a_{2})}+2^{1-(a_{1}+a_{2}+a_{3})}-2^{1-(a_{1}+a_{2}+a_{3}+a_{4})}$.
Therefore, acting by induction, we get exactly the expansion of $?(x)$. This shows that in fact $v(x)=?(x)$. $\square$\\

Now,  for $x\in(0,1]$, $x=[[b_{1},b_{2},b_{3},\ldots]]$, $\ell\geq 0$, let us introduce
\begin{eqnarray*}
\quad h_{\ell}(x)=2^{\ell-\sum_{i=1}^{\ell}b_{i}}-2^{\ell+2-\sum_{i=1}^{\ell+1}b_{i}}.
\end{eqnarray*}
(Here for $\ell=0$ an empty sum is assumed to be $0$). This is obviously a non-negative function (since $b_{\ell+1}\geq 2$). For the next identity to hold for \emph{all} real numbers in $[0,1]$, we must make a convention concerning the rational numbers. For example, if $x$ is rational, we use the finite expansion of $x$ in the definition of $h$. Then $h$ has a countable set of non-continuity points and it is continuous from the right. We immediately obtain the equality
\begin{eqnarray}
\sum\limits_{\ell=0}^{\infty}h_{\ell}(x)=1-\sum\limits_{\ell=1}^{\infty}2^{\ell-\sum_{i=1}^{\ell}b_{i}}\mathop{=}^{(\ref{sum})}1-?(x),\quad x\in(0,1].
\label{quest}
\end{eqnarray}
It appears that this series establishes exactly the relation between our formula for the moments and the semi-regular continued fractions. More precisely,
\begin{prop} For $L\geq 1$, $\ell\geq 0$, we have
\begin{eqnarray*}
L!\int\limits_{0}^{1}x^{L-1}h_{\ell}(x)\d x=\quad\int\limits_{[0,\infty)^{\ell+1}}x_{0}^{L}\cdot
\frac{(x_{0}x_{\ell})^{-1/2}\cdot\prod\limits_{i=0}^{\ell-1}I_{1}(2\sqrt{x_{i}x_{i+1}})}{\prod\limits_{i=0}^{\ell}e^{x_{i}}(2e^{x_{i}}-1)}\d x_{0}\ldots\d x_{\ell}.
\end{eqnarray*}
\end{prop}
{\it Proof.} First, given a  semi-regular continued fraction, we perform what in the theory of general continued fractions is called {\it an equivalence transformation}. In our case, we easily see that
\begin{eqnarray}
[[b_{1},b_{2},\ldots,b_{\ell+1}]]=\cfrac{b_{1}^{-1}}{1-\cfrac{(b_{1}b_{2})^{-1}}{1-\cfrac{(b_{2}b_{3})^{-1}}{\ddots 1-(b_{\ell}b_{\ell+1})^{-1}}}}:=
\Big{\langle} b_{1}^{-1},(b_{1}b_{2})^{-1},(b_{2}b_{3})^{-1},\ldots,(b_{\ell}b_{\ell+1})^{-1}\Big{\rangle}.\label{not}
\end{eqnarray}
 For $s\in\mathbb{N}$, let us introduce
$c_{s}=2\sum\limits_{n=2}^{\infty}\frac{1}{2^{n}n^{s}}=2\text{Li}_{s}(\frac{1}{2})-1$, where $\text{Li}_{\star}(\star)$ stands for the polylogarithm. We have:
\begin{eqnarray*}
c_{s}=\frac{1}{(s-1)!}\sum\limits_{n=2}^{\infty}2^{-n+1}\int\limits_{0}^{\infty}x^{s-1}e^{-nx}\d x =\frac{1}{(s-1)!}\int\limits_{0}^{\infty}\frac{x^{s-1}}{e^{x}(2e^{x}-1)}\d x,\quad c_{s}\sim 2^{-(s+1)}\text{ as }s\rightarrow\infty.
\end{eqnarray*}
\indent Now, consider the following sum
\begin{eqnarray*}
\mathcal{V}_{\ell}=\sum\limits_{q_{1},q_{2},\ldots,q_{\ell}=1}^{\infty}c_{L+q_{1}}\cdot c_{q_{1}+q_{2}}\cdots c_{q_{\ell-1}+q_{\ell}}\cdot c_{q_{\ell}}
\binom{L+q_{1}-1}{q_{1}}\binom{q_{1}+q_{2}-1}{q_{2}}\cdots\binom{q_{\ell-1}+q_{\ell}-1}{q_{\ell}}.
\end{eqnarray*}
(In case $\ell=0$, this means $\mathcal{V}_{0}=c_{L}$).\\
\indent {\it Remark}. Here and in the sequel, the infinite series and integrals involve only positive numbers and functions, and this is a big advantage of our approach. The convergence of these is always immediate. Soon we will see that $\sum_{\ell=0}^{\infty}\mathcal{V}_{\ell}=m_{L}$. As can be guessed, this series was not invented out of nowhere but rather was derived from the other results by the author; in particular, we already had an infinite system of linear equations which the moments do satisfy; see \cite{alkauskas1}, Proposition 5. Unfortunately, the system presented there has alternating terms and is not suitable for iteration. Nevertheless, if we start from the functional equation $(12)$  rather than from $(14)$ as is given in \cite{alkauskas1} (the fact which was, sadly, overseen), we get another system with strictly positive terms, and its iteration produces exactly the series above. Here we present an independent exposition from the scratch, and all we have to do is to make the ends meet.\\
\indent Thus, on the one hand, we have
\begin{eqnarray*}
\mathcal{V}_{\ell}&=&\int\limits_{0}^{\infty}\int\limits_{0}^{\infty}\ldots\int\limits_{0}^{\infty}
\frac{\d x_{0}}{x_{0}e^{x_{0}}(2e^{x_{0}}-1)}\frac{\d x_{1}}{x_{1}e^{x_{1}}(2e^{x_{1}}-1)}\cdots\frac{\d x_{\ell}}{x_{\ell}e^{x_{\ell}}(2e^{x_{\ell}}-1)}\\
&\times&\sum\limits_{q_{1},q_{2},\ldots,q_{\ell}=1}^{\infty}x_{0}^{L+q_{1}}x_{1}^{q_{1}+q_{2}}\cdots x_{\ell-1}^{q_{\ell-1}+q_{\ell}}x_{\ell}^{q_{\ell}}\\
&\times&\frac{1}{q_{1}!(L-1)!}\frac{1}{q_{2}!(q_{1}-1)!}\cdots\frac{1}{q_{\ell}!(q_{\ell-1}-1)!}\frac{1}{(q_{\ell}-1)!}.
\end{eqnarray*}
So, this $\ell-$fold sum splits. Note that the sum
\begin{eqnarray*}
\sum\limits_{q=1}^{\infty}\frac{x^q}{(q-1)!q!}=\sqrt{x}I_{1}(2\sqrt{x}),
\end{eqnarray*}
where $I_{1}(x)=J_{1}(ix)$ is the Bessel $I-$function \cite{watson}. Therefore, $(L-1)!\mathcal{V}_{\ell}$ equals precisely the integral in the formulation of the proposition.\\

Further, let $d_{1}=b_{1}^{-1}$, $d_{2}=(b_{1}b_{2})^{-1}$,$\ldots$, $d_{\ell+1}=(b_{\ell}b_{\ell+1})^{-1}$ for simplicity. Now, let us unfold the series for $\mathcal{V}_{\ell}$. This gives:
\begin{eqnarray*}
\mathcal{V}_{\ell}&=&\sum\limits_{q_{1},q_{2},\ldots,q_{\ell}=1}^{\infty}\quad\sum\limits_{b_{1},b_{2},\ldots,b_{\ell+1}=2}^{\infty}2^{\ell+1-(b_{1}+b_{2}+\ldots+b_{\ell+1})}
\frac{1}{b_{1}^{L+q_{1}}}\cdot\frac{1}{b_{2}^{q_{1}+q_{2}}}\cdots\frac{1}{b_{\ell}^{q_{\ell-1}+q_{\ell}}}\cdot\frac{1}{b_{\ell+1}^{q_{\ell}}}\\
&\times&\binom{L+q_{1}-1}{q_{1}}\binom{q_{1}+q_{2}-1}{q_{2}}\cdots\binom{q_{\ell-1}+q_{\ell}-1}{q_{\ell}}\\
&=&\sum\limits_{b_{1},b_{2},\ldots,b_{\ell+1}=2}^{\infty}2^{\ell+1-(b_{1}+b_{2}+\ldots+b_{\ell+1})}\sum\limits_{q_{1},q_{2},\ldots,q_{\ell}=1}^{\infty}
d_{1}^{L}\cdot d_{2}^{q_{1}}\cdots d_{\ell+1}^{q_{\ell}}\\
&\times&\binom{L+q_{1}-1}{q_{1}}\binom{q_{1}+q_{2}-1}{q_{2}}\cdots\binom{q_{\ell-1}+q_{\ell}-1}{q_{\ell}}.
\end{eqnarray*}
Let us perform the summation with respect to $q_{\ell}$. Note the binomial formula
\begin{eqnarray*}
\Big{(}\frac{1}{1-x}\Big{)}^{N}=\sum\limits_{s=0}^{\infty}\binom{s+N-1}{s}x^{s}\text{ for }
N\in\mathbb{N}, \quad |x|< 1.
\end{eqnarray*}
Thus, according to this expansion, we have
\begin{eqnarray*}
\sum\limits_{q_{\ell}=1}^{\infty}\binom{q_{\ell-1}+q_{\ell}-1}{q_{\ell}}d_{\ell+1}^{q_{\ell}}=
\Big{(}\frac{1}{1-d_{\ell+1}}\Big{)}^{q_{\ell-1}}-1=\langle1,d_{\ell+1}\rangle^{q_{\ell-1}}-\langle1\rangle^{q_{\ell-1}},
\end{eqnarray*}
where $\langle\star\rangle$ is the notation introduced by (\ref{not}).
Further, the summation  with respect to $q_{\ell-1}$ gives
\begin{eqnarray*}
\sum\limits_{q_{\ell-1}=1}^{\infty}\binom{q_{\ell-2}+q_{\ell-1}-1}{q_{\ell-1}}\Big{[}\Big{(}\frac{d_{\ell}}{1-d_{\ell+1}}\Big{)}^{q_{\ell-1}}-
d_{\ell}^{q_{\ell-1}}\Big{]}=\langle1,d_{\ell},d_{\ell+1}\rangle^{q_{\ell-2}}-\langle1,d_{\ell}\rangle^{q_{\ell-2}}.
\end{eqnarray*}
Now it is obvious that this procedure can be iterated. Eventually, this implies
\begin{eqnarray}
\mathcal{V}_{\ell}&=&\sum\limits_{b_{1},b_{2},\ldots,b_{\ell+1}=2}^{\infty}2^{\ell+1-(b_{1}+b_{2}+\ldots+b_{\ell+1})}
\Big{(}\langle d_{1},d_{2},\ldots,d_{\ell+1}\rangle^{L}-\langle d_{1},d_{2},\ldots,d_{\ell}\rangle^{L}\Big{)}\nonumber\\
&=&\sum\limits_{b_{1},b_{2},\ldots,b_{\ell+1}=2}^{\infty}2^{\ell+1-(b_{1}+b_{2}+\ldots+b_{\ell+1})}
\Big{(}[[b_{1},b_{2},\ldots,b_{\ell+1}]]^{L}-[[b_{1},b_{2},\ldots,b_{\ell}]]^{L}\Big{)}\nonumber\\
&=&\mathcal{A}_{\ell+1}-\mathcal{A}_{\ell}>0.\label{suma}
\end{eqnarray}
(Note that the summation with respect to $b_{\ell+1}$ was performed for the second summand). Here for brevity we put
\begin{eqnarray}
\mathcal{A}_{\ell}=\sum\limits_{b_{1},b_{2},\ldots,b_{\ell}=2}^{\infty}2^{\ell-(b_{1}+b_{2}+\ldots+b_{\ell})}
[[b_{1},b_{2},\ldots,b_{\ell}]]^{L},\text{ for }\ell\geq 1, \quad \mathcal{A}_{0}=0.\label{tarp}
\end{eqnarray}
\\

Finally, consider the function $f_{\ell+1}([[b_{1},b_{2},\ldots,b_{\ell+1},\ldots]])=2^{\ell+1-(b_{1}+b_{2}+\ldots+b_{\ell+1})}$ (for rational argument we make the same convention as in the definition of $h$; thus, we use finite expansion). This function is constant in the interval $(x,y)$, where $x=[[b_{1},b_{2},\ldots,b_{\ell+1}]]$, $y=[[b_{1},b_{2},\ldots,b_{\ell+1}-1]]$, and takes the value $f_{\ell+1}(x)$.  Therefore,
\begin{eqnarray*}
L\int\limits_{0}^{1}f_{\ell+1}(x)x^{L-1}\d x&=&
L\sum\limits_{b_{1},b_{2},\ldots,b_{\ell+1}=2}^{\infty}2^{\ell+1-(b_{1}+b_{2}+\ldots+b_{\ell+1})}
\int\limits_{[[b_{1},b_{2},\ldots,b_{\ell+1}]]}^{[[b_{1},b_{2},\ldots,b_{\ell+1}-1]]}x^{L-1}\d x\\
&=&\sum\limits_{b_{1},b_{2},\ldots,b_{\ell+1}=2}^{\infty}2^{\ell+1-(b_{1}+b_{2}+\ldots +b_{\ell+1})}\\
&\times&\Big{(}[[b_{1},b_{2},\ldots,b_{\ell+1}-1]]^{L}-[[b_{1},b_{2},\ldots,b_{\ell+1}]]^{L}\Big{)}\\
&=&\sum\limits_{b_{1},b_{2},\ldots,b_{\ell}=2}^{\infty}2^{\ell-1-(b_{1}+b_{2}+\ldots +b_{\ell})}[[b_{1},b_{2},\ldots,b_{\ell},1]]^{L}-\frac{\mathcal{A}_{\ell+1}}{2}.
\end{eqnarray*}
Each number $[[b_{1},b_{2},\ldots,b_{\ell},1]]$ can be written either in the form
$[[b_{1},b_{2},\ldots,b_{\ell-i},\underbrace{2,2,\ldots,2}_{i},1]]$ = $[[b_{1},b_{2},\ldots,b_{\ell-i}-1]]$, $0\leq i\leq\ell-1$, $b_{\ell-i}\geq 3$, or it is $[[\underbrace{2,2,\ldots,2}_{\ell},1]]=1$. Thus, this implies
\begin{eqnarray}
L\int\limits_{0}^{1}f_{\ell+1}(x)x^{L-1}\d x=-\frac{\mathcal{A}_{\ell+1}}{2}+
\sum\limits_{i=0}^{\ell-1}\frac{\mathcal{A}_{\ell-i}}{2^{i+2}}+\frac{1}{2^{\ell+1}}\mathop{<}^{(\ref{suma})}\frac{1}{2^{\ell+1}}.\label{ineq}
\end{eqnarray}
Consequently,
\begin{eqnarray*}
L\int\limits_{0}^{1}\Big{[}f_{\ell}(x)-2f_{\ell+1}(x)\Big{]}x^{L-1}\d x=\mathcal{A}_{\ell+1}-\mathcal{A}_{\ell}=\mathcal{V}_{\ell},\quad \ell\geq 0.
\end{eqnarray*}
Since $f_{\ell}(x)-2f_{\ell+1}(x)=h_{\ell}(x)$, this establishes the validity of the proposition. Moreover, 
\begin{eqnarray*}
L\int\limits_{0}^{1}h_{\ell}(x)x^{L-1}\d x\leq L\int\limits_{0}^{1}f_{\ell}(x)x^{L-1}\d x\mathop{<}^{(\ref{ineq})}\frac{1}{2^{\ell}}.\quad\square
\end{eqnarray*}

{\it Proof of the Theorem.} We see that
\begin{eqnarray*}
L!\int\limits_{0}^{1}x^{L-1}\sum\limits_{\ell=0}^{\infty}h_{\ell}(x)\d x
\mathop{=}^{(\ref{quest})}L!\int\limits_{0}^{1}x^{L-1}[1-?(x)]\d x=(L-1)!\int\limits_{0}^{1}x^{L}\d ?(x).\quad\square
\end{eqnarray*}
\indent \emph{A posteriori}, since $\sum\limits_{i=0}^{\ell}\mathcal{V}_{i}=\mathcal{A}_{\ell+1}$ as is given by (\ref{tarp}), the moments can also be expressed by the limit $m_{L}=\lim_{\ell\rightarrow\infty}\mathcal{A}_{\ell}$, and this formula is in many respects very different from (\ref{riba}).

\par\bigskip

\noindent {\sc Giedrius Alkauskas}, Institute of Mathematics, Department of Integrative Biology,
Universit\"{a}t f\"{u}r Bodenkultur Wien, Gregor Mendel-Stra{\ss}e 33, A-1180 Wien, Austria
\\
{\tt giedrius.alkauskas@gmail.com}


\begin{thebibliography}{9}
\bibitem{alkauskas1}G. Alkauskas, {\it The moments of Minkowski question mark function: the dyadic period function},
 Glasg. Math. J. {\bf 52} (1) (2010), 41--64.

\bibitem{alkauskas2}G. Alkauskas, {\it Generating and zeta functions, structure, spectral and analytic properties of the moments
of the Minkowski question mark function}, Involve {\bf 2} (2) (2009), 121--159.

\bibitem{alkauskas3}G. Alkauskas, {\it An asymptotic formula for the moments of Minkowski question mark function in
the interval $[0,1]$}, Lith. Math. J. {\bf 48} (4) (2008), 357--367.

\bibitem{alkauskas4}G. Alkauskas, {\it The Minkowski question mark function: explicit series for the dyadic period function and moments},
Math. Comp. {\bf 79} (269) (2010), 383--418.

\bibitem{grabner_tichy}P. J. Grabner, P. Kirschenhofer and R. F. Tichy,
{\it Combinatorial and arithmetical properties of linear numeration systems}, Combinatorica {\sc 22} (2) (2002), 245--267.

\bibitem{kessebohmer}M. Kesseb\"{o}hmer, B. O. Stratmann,
{\it Fractal analysis for sets of non-differentiability of Minkowski's
question mark function}, J. Number Theory {\bf 128} (9) (2008), 2663-2686.

\bibitem{lamberger}M. Lamberger, {\it On a family of singular measures related
to Minkowski's $?(x)$ function}, Indag. Mathem., N.S., {\bf 17} (1) (2006), 45--63.

\bibitem{panti} G. Panti, {\it Multidimensional continued fractions and a
    Minkowski function}, Monatsh. Math., {\bf 154} (3) (2008), 247-264.

\bibitem{ramharter}G. Ramharter, {\it On Minkowski's singular function}, Proc. Amer. Math. Soc.  {\bf 99} (3) (1987), 596--597.

\bibitem{ryde} F. Ryde, {\it On the relation between two Minkowski functions}, J. Number Theory,  {\bf 17} (1) (1983), 47--51.

\bibitem{watson}G. N. Watson, {\it A treatise on the theory of Bessel functions, 2nd ed.} Cambridge University Press, 1996.

\bibitem{mink} An extensive bibliography on the Minkowski question mark function, available at
\url{http://www.alkauskas.puslapiai.lt/minkowski.htm}

\end{thebibliography}
\end{document}